\documentclass[11pt,fleqn]{article}
\usepackage{amsfonts}
\newcommand{\ol}{\setlength{\itemsep}{0pt.}\begin{enumerate}}
\newcommand{\eol}{\end{enumerate}\setlength{\itemsep}{-\parsep}}
\newcommand{\ignore}[1]{}
\title{A modified logarithmic Sobolev inequality for the Hamming cube
  and some applications}
\author{Alex Samorodnitsky}

\begin{document}
\date{}
\maketitle


\newtheorem{THEOREM}{Theorem}[section]
\newenvironment{theorem}{\begin{THEOREM} \hspace{-.85em} {\bf :}
}%
                        {\end{THEOREM}}
\newtheorem{LEMMA}[THEOREM]{Lemma}
\newenvironment{lemma}{\begin{LEMMA} \hspace{-.85em} {\bf :} }%
                      {\end{LEMMA}}
\newtheorem{COROLLARY}[THEOREM]{Corollary}
\newenvironment{corollary}{\begin{COROLLARY} \hspace{-.85em} {\bf
:} }%
                          {\end{COROLLARY}}
\newtheorem{PROPOSITION}[THEOREM]{Proposition}
\newenvironment{proposition}{\begin{PROPOSITION} \hspace{-.85em}
{\bf :} }%
                            {\end{PROPOSITION}}
\newtheorem{DEFINITION}[THEOREM]{Definition}
\newenvironment{definition}{\begin{DEFINITION} \hspace{-.85em} {\bf
:} \rm}%
                            {\end{DEFINITION}}
\newtheorem{EXAMPLE}[THEOREM]{Example}
\newenvironment{example}{\begin{EXAMPLE} \hspace{-.85em} {\bf :}
\rm}%
                            {\end{EXAMPLE}}
\newtheorem{CONJECTURE}[THEOREM]{Conjecture}
\newenvironment{conjecture}{\begin{CONJECTURE} \hspace{-.85em}
{\bf :} \rm}%
                            {\end{CONJECTURE}}
\newtheorem{MAINCONJECTURE}[THEOREM]{Main Conjecture}
\newenvironment{mainconjecture}{\begin{MAINCONJECTURE} \hspace{-.85em}
{\bf :} \rm}%
                            {\end{MAINCONJECTURE}}
\newtheorem{PROBLEM}[THEOREM]{Problem}
\newenvironment{problem}{\begin{PROBLEM} \hspace{-.85em} {\bf :}
\rm}%
                            {\end{PROBLEM}}
\newtheorem{QUESTION}[THEOREM]{Question}
\newenvironment{question}{\begin{QUESTION} \hspace{-.85em} {\bf :}
\rm}%
                            {\end{QUESTION}}
\newtheorem{REMARK}[THEOREM]{Remark}
\newenvironment{remark}{\begin{REMARK} \hspace{-.85em} {\bf :}
\rm}%
                            {\end{REMARK}}

\newcommand{\thm}{\begin{theorem}}
\newcommand{\lem}{\begin{lemma}}
\newcommand{\pro}{\begin{proposition}}
\newcommand{\dfn}{\begin{definition}}
\newcommand{\rem}{\begin{remark}}
\newcommand{\xam}{\begin{example}}
\newcommand{\cnj}{\begin{conjecture}}
\newcommand{\mcnj}{\begin{mainconjecture}}
\newcommand{\prb}{\begin{problem}}
\newcommand{\que}{\begin{question}}
\newcommand{\cor}{\begin{corollary}}
\newcommand{\prf}{\noindent{\bf Proof:} }
\newcommand{\ethm}{\end{theorem}}
\newcommand{\elem}{\end{lemma}}
\newcommand{\epro}{\end{proposition}}
\newcommand{\edfn}{\bbox\end{definition}}
\newcommand{\erem}{\bbox\end{remark}}
\newcommand{\exam}{\bbox\end{example}}
\newcommand{\ecnj}{\bbox\end{conjecture}}
\newcommand{\emcnj}{\bbox\end{mainconjecture}}
\newcommand{\eprb}{\bbox\end{problem}}
\newcommand{\eque}{\bbox\end{question}}
\newcommand{\ecor}{\end{corollary}}
\newcommand{\eprf}{\bbox}
\newcommand{\beqn}{\begin{equation}}
\newcommand{\eeqn}{\end{equation}}
\newcommand{\wbox}{\mbox{$\sqcap$\llap{$\sqcup$}}}
\newcommand{\bbox}{\vrule height7pt width4pt depth1pt}
\newcommand{\qed}{\bbox}
\def\sup{^}

\def\H{\{0,1\}^n}
\def\B{\{0,1\}}

\def\S{S(n,w)}

\def\n{\lfloor \frac n2 \rfloor}

\def\Tp{Tchebyshef polynomial}
\def\Tps{TchebysDeto be the maximafine $A(n,d)$ l size of a code with distance $d$hef polynomials}
\newcommand{\rarrow}{\rightarrow}

\newcommand{\larrow}{\leftarrow}

\overfullrule=0pt
\def\setof#1{\lbrace #1 \rbrace}

\def \E{\mathbb E}
\def \R{\mathbb R}
\def \Z{\mathbb Z}
\def \F{\mathbb F}

\def\<{\left<}
\def\>{\right>}
\def \({\left(}
\def \){\right)}
\def \e{\epsilon}
\def \d{\delta}
\def \l{\lambda}

\def \z{\exp\left\{-\Omega(n)\right\}}
\def \k{O\(\exp\left\{-KN\right\}\)}
\def \A{A^{(r)}}
\def \P{{\cal P}}
\def \Q{{\cal Q}}
\def \M{{\cal M}}

\begin{abstract}
The logarithmic Sobolev inequality \cite{Gr} for the Hamming cube $\{0,1\}^n$ states that for any real-valued function $f$ on the cube holds $$\E_x \sum_{y \sim x} (f(x) - f(y))^2 \ge 2 \cdot Ent\(f^2\)$$ We show that the constant $C = 2$ at the right hand side of this inequality can be replaced by a function $C(\rho)$ depending on $\rho = \frac1n \frac{Ent\(f^2\)}{\E f^2}$. The function $C(\cdot)$ is an increasing convex function taking $[0,\log 2]$ to $[2, 2/\log 2]$.

We present some applications of this modified inequality. In particular, it is used to obtain a discrete version of the Faber-Krahn inequality for small subsets of the Hamming cube, answering a question of Friedman and Tillich \cite{FT}. We introduce, following \cite{FT}, the notion of a fractional edge-boundary size of a subset of $\H$, and show Hamming balls of radius at most $n/2 - \tilde{\Omega}\(n^{3/4}\)$ to be sets with (asymptotically) the smallest fractional edge-boundary for their size.
\end{abstract}

\section{Introduction}

\subsection{Isoperimetric problems on the Hamming cube}

This paper deals with {\it discrete isoperimetric inequalities} on graphs. Let a graph $G = (V,E)$ be given, and let $A \subseteq V$ be a set of vertices in $G$. An isoperimetric inequality addresses the question of how small the boundary $\partial A$ of $A$ can be, given the cardinality of $A$, by lower bounding the size of $\partial A$ by an appropriate function of $|A|$.  There are various ways do define and measure the boundary of the set, two salient examples being the {\it vertex boundary} of $A$, consisting of
the vertices of $A$ which have neighbors outside $A$, and the {\it edge
boundary} of $A$, which is the set of edges crossing from $A$ to its
complement. In both these cases, the size of the boundary is the cardinality of the corresponding set of vertices (or edges). However, one can also think of examples of a somewhat different nature, some of which will be considered below.

In this paper we deal with
a specific graph - the Hamming cube $\H$. This is a graph with $2^n$ vertices
indexed by boolean strings of length $n$. Two vertices are
connected by an edge if they differ only in one coordinate. The
metric defined by this graph is called the Hamming distance. In other words, two vertices $x$
and $y$ are at distance $d$ if they differ in $d$ coordinates. Let us define two important families of subsets of $\H$. A {\it Hamming ball} is a ball in the Hamming metric. A {\it subcube} is a
subset of the vertices obtained by fixing the value in some of the
coordinates. The number of fixed coordinates is called the co-dimension of
the subcube. It turns out \cite{Harper,Harp, Hart} that a Hamming ball has the smallest vertex boundary for its
size, and a subcube has the smallest edge boundary.

We will consider several versions of the edge-isoperimetric inequality
for the cube. Let $|\partial A|$ be the cardinality of the edge
boundary of $A$, normalized, for convenience, by $2^{n-1}$.
The standard version of the inequality \cite{Harp, Hart} states that for
any subset $A \subseteq \H$ we have
\footnote{We use natural logarithms throughout the paper.}
\beqn
\label{set_isop}
|\partial A|~ \ge~ \frac{2}{\log 2} \cdot \frac{|A|}{2^n} \log
\frac{2^n}{|A|}.
\eeqn
This is tight if $A$ is a subcube of an arbitrary co-dimension $0\le t
\le n$.

\noindent {\it A logarithmic Sobolev inequality} \cite {Gr} establishes a relation
between two (appropriately defined) notions: the variation of a function
and its entropy. For a function $f$ on the Hamming cube $\H$ endowed
with the uniform measure, this
translates to
\beqn
\label{ineq_Sob}
\E_x \sum_{y \sim x} (f(x) - f(y))^2 \ge 2 \cdot Ent\(f^2\) = 2\cdot \(\E f^2
\log f^2 - \E f^2 \log \E f^2\)
\eeqn
It is useful to view (\ref{ineq_Sob}) as an isoperimetric
inequality. In particular (\cite{FS}), it implies a functional form of the
edge-isoperimetric inequality (\ref{set_isop}). For a non-zero
function $f:~\H \rarrow \R$ holds
\beqn
\label{func_isop}
\E_x \sum_{y \sim x} (f(x) - f(y))^2 \ge 2 \cdot \E f^2 \log \frac{\E
  f^2}{\E^2 |f|}
\eeqn
Choosing $f$ in (\ref{func_isop}) to be the
characteristic function of a subset $A$ of the cube, we recover the
edge-isoperimetric inequality (\ref{set_isop}) with a somewhat worse
constant on the right hand side, replacing $2/\log 2$ with $2$. On the
other hand, there are examples (\cite{FS}) of real-valued functions $f$ for which the constant $2$ in (\ref{func_isop}) is tight.

Let us now briefly sketch the results in this paper, before discussing them in fuller detail in subsection~\ref{subsec:res} below. We will show that the constant $C = 2$ at the right hand side of the logarithmic Sobolev inequality (\ref{ineq_Sob}) can be replaced by a function $C(\rho)$ depending on $\rho = \frac1n \frac{Ent\(f^2\)}{\E f^2}$. The function $C(\rho)$ will be given explicitly. It is a convex function which increases from $2$ to $2/\log 2$ as $\rho$ goes from $0$ to $\log 2$. We will also observe that $C(\rho)$ gives the correct dependence on $\rho$, describing functions on $\{0,1\}^n$ for which the modified inequality is tight.

This will imply a corresponding modification of the functional isoperimetric inequality (\ref{func_isop}). Here, as well as in (\ref{ineq_Sob}), it will be possible to replace the constant $C = 2$ by the function $C(\rho)$.

The modified version of (\ref{func_isop}) is used to derive a discrete analogue of the classical Faber-Krahn inequality
in $\R^n$ for small subsets of the Hamming cube $\H$, answering\footnote{Up to an error which becomes negligible as the dimension of the cube grows.} a question from \cite{FT}. We will introduce, following \cite{FT}, the notion of a fractional edge-boundary size of a subset of $\H$, and show Hamming balls of radius at most $n/2 - \tilde{\Omega}\(n^{3/4}\)$ to be sets with (asymptotically) the smallest fractional edge-boundary for their size.

The question in \cite{FT} is a part of an approach to obtain upper bounds on the cardinality of binary error-correcting codes. We will now take a brief detour to the theory of error-correcting codes in order to provide a natural framework in which this question can be discussed.

\subsection{Bounds on binary error correcting codes}

A binary error-correcting code of length $n$ and minimal distance $d$
is a subset $C$ of the boolean cube $\H$ such that the distance between
any two distinct points in $C$ is at least $d$. In other words, the
points in $C$ can be taken as centers in a disjoint packing of Hamming
balls of raduis $\lceil \frac{d-1}{2} \rceil$ into $\H$.

The question of the maximal possible cardinality $A(n,d)$ of such a packing is
one of the central questions of coding theory. The best known upper
bounds on $A(n,d)$ were obtained in \cite{MRRW} following Delsarte's
linear programming approach \cite{dels}. The analysis in \cite{MRRW}
uses theory of orthogonal polynomials and is somewhat complicated.

A different approach to obtain some of the bounds in \cite{MRRW}
was presented in \cite{FT}. The appeal of
this new approach is in showing the possibility to work with Delsarte's linear
inequalities without resorting to language and tools of orthogonal
polynomial theory. In particular, \cite{FT} establishes a connection
between packing bounds and isoperimetric questions in the Hamming
cube.
To describe this connection, we need a notion of the {\it fractional
edge boundary size} of a subset of
the cube \footnote{This is a reformulation of a closely related notion of the maximal eigenvalue of
a subset introduced in \cite{FT}.}.

For $A \subseteq \H$, the fractional edge boundary size of $A$ is defined as
\beqn
\label{dfn:frac-bound}
|\partial^* A| = \min\left\{\E_x \sum_{y\sim x} (f(x)-f(y))^2:~\mbox{supp}(f) \subseteq
A,~\E f^2 = \frac{|A|}{2^n}\right\}
\eeqn
The right hand side of this definition computes a minimum of the variation of $f$
over a certain family of functions. Note that this family contains the characteristic
function of $A$ whose variation equals the (normalized) cardinality of the edge boundary $\partial A$.
Consequently, the fractional boundary of $A$ is at most as large as its boundary.

The following result is proved in \cite{NS, NS1}, following the approach of \cite{FT}. Let us mention that this claim was proved in \cite{FT} for an important special case of linear codes.

\thm
\label{thm:ns}
Let $A$ be a subset of the boolean cube $\H$ such that
$$
|\partial^* A| \le (2d+1) \cdot \frac{|A|}{2^{n-1}}
$$
Let $C$ be a binary error-correcting code with minimal distance $d$. Then
$$
|C| \le n |A|
$$
\ethm

This suggests a way to obtain upper bounds for codes by finding
subsets of the cube with a small fractional boundary. Natural
candidates to try are the isoperimetric sets, that is Hamming balls
and subcubes. Their fractional boundaries are analyzed in
\cite{FT}. It turns out that among these two options, a Hamming ball has the smaller fractional
edge boundary. Note that this pinpoints an intriguing difference between
the notions of the fractional edge-boundary and that of the 'ordinary' edge-boundary, for which the subcubes are the optimal sets.

Let $B$ denote a Hamming ball of radius $r$. It is shown in \cite{FT} that
\beqn
\label{ineq:ball_bound}
|\partial^* B | \le 4\(\frac n2 - \sqrt{r(n-r)} +o(n)\) \cdot \frac{|B|}{2^n}
\eeqn
Combined with Theorem~\ref{thm:ns}, this shows that a
binary error-correcting code with minimal distance $d$ is at most as large, up to
negligible multiplicative factors, as a Hamming ball of radius $r =
n/2 - \sqrt{d(n-d)}$. This provides an alternative proof of the {\it first
linear programming bound} for binary codes \cite{MRRW}.

This concludes our detour into coding theory. We are now ready to state the isoperimetric problem of \cite{FT}.

\subsection{An isoperimetric problem for the Hamming cube}
In order to obtain the best possible bounds on codes via Theorem~\ref{thm:ns}, we need to find subsets of the Hamming cube with the smallest possible fractional edge-boundary. In particular, an existence of subsets whose fractional boundary is noticeably smaller than that of Hamming balls of the same cardinality, would imply an improvement on the best currently known bounds. This naturally leads to the following questions \cite{FT}

{\bf A fractional edge-isoperimetric problem for the Hamming cube}:
\begin{itemize}
\item
What is the smallest possible fractional
boundary of a subset of $\H$ of a given cardinality?
\item
Which sets have the smallest fractional boundaries?
\end{itemize}

These questions were the starting point of our investigation. Before describing our results, let us mention a connection to the classical Faber-Krahn inequality in $\R^n$, as pointed out in \cite{FT}.

First, here is a brief description of the Euclidean space inequality, following \cite{Chavel}.
For an open set $\Omega$ in $\R^n$, consider the functional
\beqn
F[\phi] = \frac{\|grad ~\phi\|^2_2}{\|\phi\|^2_2}
\label{cont_fk}
\eeqn
where $\phi$ ranges over smooth functions supported in $\Omega$, and
the associated infimum
$
\l^*(\Omega) = inf_{\phi} F[\phi]
$.

$\l^*(\Omega)$ is referred to as the fundamental tone of
$\Omega$. The Faber-Krahn inequality states that among all sets
$\Omega$ of the same measure, Euclidean ball has the minimal
fundamental tone.

In the discrete setting of the Hamming cube, a reasonable
interpretation of (\ref{cont_fk}) is to consider the functional
$$
F[f] = \frac{\E_x \sum_{y\sim x} (f(x) - f(y))^2}{\E f^2}
$$
where $f$ ranges over functions supported in a subset $A$ of $\H$. In
our terminology, the ``fundamental tone'' of $A$ is
$$
\lambda^*(A) = \min_f F[f] = \frac{2^n}{|A|} \cdot |\partial^* A|
$$
Hence, the set with the smallest fractional boundary for its size has the minimal fundamental tone, and vice versa.

Following \cite{FT} we will refer to the fractional edge-isoperimetric
problem as the discrete Faber-Krahn problem for the Hamming cube.

\subsection{Main results}
\label{subsec:res}
Our main technical result is a modified version of the logarithmic Sobolev inequality (\ref{ineq_Sob}). Let $H(x) = -x\log x - (1-x) \log(1-x)$ be the "natural" (i.e., using natural logarithms) entropy function.
\thm
\label{thm:log-Sobolev}
\begin{itemize}
\item
Let $f:~\H \rarrow \R$ be a non-zero function, and let $\rho = \frac 1n \frac{Ent\(f^2\)}{\E f^2}$. Then
\beqn
\label{ineq:l_sob}
\E_x \sum_{y \sim x} (f(x) - f(y))^2 \ge C(\rho) \cdot Ent\(f^2\),
\eeqn
where
$$
C(x) = \frac{4}{x} \cdot \(\frac12 - \sqrt{H^{-1}(\log 2 - x)\Big (1 -
  H^{-1}(\log 2 - x)\Big )}\)
$$
\item
The function $C(\cdot)$ is an increasing convex function, taking $[0, \log 2]$ to
$[2,2/\log 2]$.
\end{itemize}
\ethm

Inequality (\ref{ineq:l_sob}) is tight in
the following sense: for each $\rho \in [0, \log 2]$ there exists a non-constant
function $f = f_{\rho}$ such that $Ent\(f^2\) \ge \rho n \E f^2$ and
\beqn
\label{ineq:sob-tight}
\E_x \sum_{y \sim x} (f(x) - f(y))^2 \le (1 + o_n(1)) \cdot C(\rho)
\cdot Ent\(f^2\)
\eeqn
This follows from the tightness of inequality (\ref{ineq:fk}) below (see the second part of Theorem~\ref{thm:fk}), since that inequality is a corollary of (\ref{ineq:l_sob}). The fact that (\ref{ineq:fk}) is tight for Hamming balls follows from (\ref{ineq:ball_bound}). The functions $f_{\rho}$ are the minimal variation functions supported on a Hamming ball of an appropriate radius. They are constructed explicitly in \cite{FT}.

Theorem~\ref{thm:log-Sobolev} together with the observation $Ent\(f^2\) \ge \E f^2 \log \frac{\E f^2}{\E^2 |f|}$ (\cite{FS}), implies a corresponding modification of the functional isoperimetric inequality (\ref{func_isop}).
\cor
\label{cor:ineq:isop}
\beqn
\label{ineq:isop}
\E_x \sum_{y \sim x} (f(x) - f(y))^2 \ge C(\rho) \cdot \E f^2 \log,
\frac{\E f^2}{\E^2 |f|}
\eeqn
\ecor
where $\rho = \frac 1n \log \frac{\E f^2}{\E^2 |f|}$. Inequality (\ref{ineq:isop}) is tight in the same sense and for the same reasons (\ref{ineq:sob-tight}) is tight.

Let us briefly discuss this inequality. It implies, in particular, that as the ratio $\frac{\E f^2}{\E^2 f}$ grows (the function $f$ becomes less "flat") its edge-isoperimetric constant approaches the isoperimetric constant $C = \frac{2}{\log 2}$ in the edge-isoperimetric inequality (\ref{func_isop}) for $0$-$1$ functions. One possible partial explanation for this phenomenon is that, for functions supported on a small set $A \subseteq \H$, the main contribution to the variation $\E_x \sum_{y\sim x} (f(x) - f(y))^2$ is likely to come from edges $(x,y)$ which belong to the edge-boundary of $A$.

It is now straightforward to derive the fractional edge-isoperimetric inequality (\ref{ineq:fk}) from Corollary~\ref{cor:ineq:isop}.
Let $f$ be a function supported on a subset $A$ of $\H$. Then, by the Cauchy-Schwarz inequality, $\frac{\E f^2}{\E^2 f} \ge \frac{2^n}{|A|}$. Since the function $C(\cdot)$ is monotone, we have
$$
\E_x \sum_{y\sim x} (f(x)-f(y))^2 \ge C\(\frac 1n \log \frac{2^n}{|A|}\) \cdot \log \frac{2^n}{|A|} \cdot \E f^2
$$
Recalling the definition of the fractional edge boundary (\ref{dfn:frac-bound}), and substituting the explicit expression for $C(\cdot)$, gives (\ref{ineq:fk}).

The inequality (\ref{ineq:fk}), together with the second part of Theorem~\ref{thm:fk},
provide an asymptotic solution to the Faber-Krahn problem for the Hamming cube, at least in the range of interest to the coding theory. It turns out that, up to an error which becomes negligible as the dimension $n$ grows, Hamming balls of radius $0 \le r \le \frac n2 -  o(n)$ are the sets with the smallest fractional boundary (fundamental tone) for their size.

From the viewpoint of coding theory, this implies that Theorem~\ref{thm:ns} cannot lead to an improvement on the best currently known bounds for binary codes \cite{MRRW}. Let us briefly discuss one implication of this fact. The bounds in \cite{MRRW} are obtained following Delsarte's linear programming approach, and there are claims in coding theory (\cite{BJ,NS}) which seem to indicate that these are the best bounds attainable with this approach. Since Theorem~\ref{thm:ns} is derived within the same linear programming framework, inequality (\ref{ineq:fk}) may be interpreted as an additional evidence in this direction\footnote{and it is not, in this sense, very surprising.}. It seems worthwhile to point out that, in this manner, the coding theory provides both the question prompting this investigation and an indication of what the answer might be, by suggesting the putative optimality of Hamming balls for the Faber-Krahn problem.

\thm
\label{thm:fk}
\begin{itemize}
\item
Let $H(x) = -x \log x - (1-x) \log(1-x)$ be the entropy function. Then for any subset $A$ of $\H$ holds
\beqn
\label{ineq:fk}
|\partial^* A| \ge 4n \(\frac12 - \sqrt{H^{-1}\(\frac{\log
    |A|}{n}\)\(1 - H^{-1}\(\frac{\log |A|}{n}\)\)}\) \cdot \frac{|A|}{2^n}
\eeqn
\item
On the other hand, let $B$ be a Hamming ball. Then
$$
|\partial^* B| \le 4n \(\frac12 - \sqrt{H^{-1}\(\frac{\log
    |B|}{n}\)\(1 - H^{-1}\(\frac{\log |B|}{n}\)\)} +o_n(1)\) \cdot \frac{|B|}{2^n}
$$
\end{itemize}
\ethm

The second part of this theorem is due to \cite{FT}. It follows from (\ref{ineq:ball_bound}) and the fact that the cardinality of a Hamming ball of radius $r$ is at least $\exp\{n H\(\frac{r}{n}\) - o(n)\}$ \cite{vLint}.

The bound in (\ref{ineq:fk}) is not very good for balanced subsets $A$, for which $\frac{\log |A|}{n}$ is close to $\log 2$. For instance, for $|A| = 2^{n-1}$, the bound gives $|\partial^* A| \ge \log 2$. On the other hand, it is not hard to see that the correct bound in this case is $|\partial^* A| \ge 1$. Equality is attained on a subcube of co-dimension one (but not on a Hamming ball of radius $n/2$).

On the other hand, (\ref{ineq:fk}) is interesting, as long as the error term $o(1)$ in the second part of Theorem~\ref{thm:fk} is of order lower than that of the main term. This error term is of order $n^{-1/2}$, up to poly-logarithmic terms\footnote{This could be derived from the computation in \cite{FT}, or, alternatively, from the estimates on the minimal roots of Krawchouk polynomials \cite{FS, NS}.}. Therefore, the theorem provides a satisfactory lower bound for the fractional edge boundary size of $A$ as long as
$$
\Big |\frac{\log |A|}{n} - \log 2 \Big | \ge \tilde{\Omega}\(n^{-1/2}\)
$$
In particular, Hamming balls of radius at most $n/2 - \tilde{\Omega}\(n^{3/4}\)$ have (asymptotically) the smallest fractional edge-boundary for their size.

\noindent {\bf Questions}:
\begin{itemize}
\item
How small can the fractional edge boundary size of a balanced subset $A$ of $\H$ be?
\item
Which balanced sets have the smallest fractional edge-boundary?
\end{itemize}

Finally, let us briefly mention an additional application of inequality (\ref{ineq:l_sob}), of a somewhat different nature. The 'standard' logarithmic Sobolev inequality (\ref{ineq_Sob}) is used in \cite{Gr} to derive a hypercontractive inequality for functions on the discrete cube (see also \cite{Bonami, Beckner}).
In order to describe this inequality, let $\{w_S\}_{S\in \H}$ be the character basis in the space of real-valued functions on the Hamming cube.
Let $0 \le t \le 1$ and let $T = T_t$ be
a linear operator taking a function $f = \sum_S \hat{f}(S) w_S$ to $Tf
= \sum_S t^{|S|} \hat{f}(S) w_S$. Then (\cite{Gr})
\beqn
\label{ineq:beckner}
\|Tf\|_2 \le \|f\|_{1+t^2}
\eeqn

Substituting  (\ref{ineq:l_sob}) instead of  (\ref{ineq_Sob}) in the proof in \cite{Gr} leads to a modified version of (\ref{ineq:beckner}). It turns out that the exponent $2$ on the right hand side of the inequality can be replaced by a function $e(\rho)$, depending on $\rho = \frac1n \frac{Ent\(f^2\)}{\E f^2}$. The function $e(\rho)$ is decreasing, with $e(0) = 2$.\footnote{Hence, for functions with high entropy, this gives a strengthening of (\ref{ineq:beckner}).} This modified inequality might be useful in coding theory, following the applications of (\ref{ineq:beckner}) in \cite{KL, ACKL}.

\section{The proof of Theorem~\ref{thm:log-Sobolev}}
Let us start with a brief overview.
The main goal of this section is to prove the logarithmic Sobolev inequality (\ref{ineq:l_sob}).
Our proof follows the outline of the proof of (\ref{ineq_Sob}) in \cite{Gr}. We will prove an inequality (\ref{ineq:tech}), which will imply (\ref{ineq:l_sob}) as a corollary, first for the base case $n=1$, and then for general $n$, using subadditivity of entropy. Compared to \cite{Gr}, we need to prove a bit more for the base case (see Remark~\ref{rem:diff} below) and to carry this additional information along to the general case. This seems to complicate things somewhat, making it necessary to go through the intermediate inequality (\ref{ineq:tech}).

We will prove the second claim of the theorem, on the properties of the function $C(\rho)$ in (\ref{ineq:l_sob}), along the way, in Lemma~\ref{lem:functions}.

We may and will assume from now on that we deal only with nonnegative functions on $\H$, since substituting $|f|$ instead of $f$ in (\ref{ineq:l_sob}) decreases the left hand side and does not affect the right hand side.

Several functions on the real line play an important role in the proof. We start by defining these functions and stating some of their properties. Let $\psi$ be defined on $[0,1]$ by
\beqn
\label{def:psi}
\psi(t) = \frac12 (1-\sqrt t)^2 \log (1-\sqrt t)^2 + \frac12
(1+\sqrt t)^2 \log (1+\sqrt t)^2 - (1+t)\log(1+t)
\eeqn
In other words, $\psi(t) = Ent\(f^2\)$, where $f$ is a function on
$\B$ with $f(0) = 1 - \sqrt t$, $f(1) = 1 + \sqrt t$.

The following main technical lemma lists the relevant properties of $\psi$ and several derived functions, including the function $C$. The proof of the lemma is rather long and is postponed till the next section.
\lem
\label{lem:functions}
\begin{enumerate}
\item
The function $\psi$ is strictly increasing and concave on $[0,1]$, taking this interval onto $[0, 2\log 2]$. This allows us to define the inverse function $\phi = \psi^{-1}$. This is a strictly increasing convex function taking $[0,2\log 2]$ onto $[0,1]$.
\item
The function $\psi(t)/(1 + t)$ is strictly increasing and concave on $[0,1]$, taking this interval onto $[0, \log 2]$. This allows us to define the inverse function $\alpha(t) = \(\frac{\psi(t)}{1 + t}\)^{-1}$. This is a strictly increasing convex function taking $[0, \log 2]$ to $[0,1]$.
\item
The function $c(t) = \frac{4 \alpha(t)}{t(1 + \alpha(t))}$ is strictly increasing and convex on $[0, \log 2]$, taking this interval onto $[2, 2/\log 2]$.\footnote{Here, as usual, we take $c(0) = \lim_{t\rarrow 0} \frac{4 \alpha(t)}{t(1 + \alpha(t))} = 2$.}
\item
The function $c(t)$ has an explicit representation
$$
c(t) = \frac{4}{t} \cdot \(\frac12 - \sqrt{H^{-1}(\log 2 - t)\Big (1 -
  H^{-1}(\log 2 - t)\Big )}\).
$$
In other words $c = C$, where $C$ is the function in (\ref{ineq:l_sob}).
\end{enumerate}
\elem
Note that the second claim of Theorem~\ref{thm:log-Sobolev} follows from the third and the fourth claims of this lemma.

\rem
\label{rem:diff}
We mentioned a difference between the proof of the base case $n=1$ here and in \cite{Gr}. Let us give some details. In \cite{Gr},  (\ref{ineq_Sob}) is shown for $\{0,1\}$, which, in our notation, amounts to proving an inequality $\psi(t) \le 2t$ on $[0,1]$. We show, in addition, that $\psi$ is concave on $[0,1]$. This additional convexity property turns out to be crucial for our proof.
\erem

Next, we introduce additional notation, and prove a simple auxiliary inequality.

For a function $f$ on $\H$ let $D^2(f) = \E_x \sum_{y \sim x} (f(x) -
f(y))^2$, and let $K^2(f) = \frac 14 \cdot \E_x \sum_{y \sim x} (f(x) +
f(y))^2$. Note that $K^2(f) = n \E f^2 - \frac14 \cdot D^2(f)$. Note
also that for a non-zero nonnegative function $f$, $K^2(f)$ is strictly positive.

\lem
\label{lem:ent_and_k}
For a nonnegative function $f$ on $\H$ holds
$$
Ent\(f^2\) \le 2\log 2 \cdot K^2(f)
$$
\elem
\prf
We will prove the claim for the base case $n=1$ and then use subadditivity of entropy to deduce it for the general case.

The case $n = 1$. We may assume $f$ is non-zero, since the claim is trivially true otherwise. Both sides of the inequality are
$2$-homogeneous, and consequently we may assume $\E f = 1$. Without loss of generality,  $f(0) = 1 - s$,
$f(1) = 1 + s$ for some $0 \le s \le 1$. We need to show
$$
Ent\(f^2\) = \psi\(s^2\) \le 2 \log 2 \cdot K^2(f) = 2\log 2,
$$
and this is indeed true by the first claim of Lemma~\ref{lem:functions}, since $\psi$ is increasing and $\psi(1) = 2\log 2$.

The general case $n \ge 1$. Let $1 \le i \le n$ be an index of a coordinate and let $x \in \H$. Fixing all the coordinates $j \not = i$ to be $x_j$ we obtain a copy of $\B$. Let $f^{(x)}_i$ be the restriction of $f$ to this one-dimensional cube.

Recall (\cite{Ledoux}) that entropy is subadditive, namely
$$
\sum_{i=1}^n \E_{x} Ent\(f^{(x)}_i\) \ge Ent(f),
$$
while $D^2(f)$ is additive, that is
$$
\sum_{i=1}^n \E_{x} D^2\(f^{(x)}_i\) = D^2(f).
$$
$K^2(f)$ is also
additive, since $\sum_{i=1}^n \E_x K^2\(f^{(x)}_i\) = \sum_{i=1}^n
\E_x \(\E \(f^{(x)}_i\)^2 - 1/4 D^2\(f^{(x)}_i\)\) = n \E f^2 - 1/4 D^2(f) = K^2(f)$.

From this and the base case, we have
$$
Ent\(f^2\) \le \sum_{i=1}^n \E_x Ent\(\(f^{(x)}_i\)^2\) \le 2\log 2 \cdot \sum_{i=1}^n \E_x K^2\(f^{(x)}_i\) = 2\log 2
\cdot K^2(f)
$$
\eprf

Let us now pass to the main technical claim in this section. Note that the right hand side in (\ref{ineq:tech}) is well-defined, due to the preceding lemma.

\pro
\label{thm:technical}
Let $f$ be a non-zero nonnegative function on $\H$. Then
\beqn
\label{ineq:tech}
D^2(f) \ge 4 \cdot K^2(f) ~\phi\(\frac{Ent\(f^2\)}{K^2(f)}\)
\eeqn
\epro
The proof follows the same outline. First we prove the base case $n = 1$.
In this case, we will show equality
$$
D^2(f) = 4 \cdot K^2(f) ~\phi\(\frac{Ent\(f^2\)}{K^2(f)}\).
$$
Indeed, since $f$ is non-zero, we may assume $\E f = 1$. This implies $K^2(f) = 1$. Thus we have to prove $D^2(f) = 4
\phi\(Ent\(f^2\)\)$. Let $0 \le s \le 1$, and $f(0) = 1 - s$, $f(1) =
1 + s$. Then $D^2(f) = 4s^2$ and
$4\phi\(Ent\(f^2\)\) = 4\phi \(\psi\(s^2\) \) = 4s^2$,
verifying the base case.

We also need to deal with a slight technicality, the case $f$ is the zero function, because in the general case below, some of the one-dimensional restrictions of $f$ might be zero. In this case, we formally define $Ent\(f^2\)/K^2(f)$ to be zero. Then (\ref{ineq:tech}) remains valid (as equality) in the one-dimensional case. It is easy to see that this formal definition does not affect the computation below.

The general case. Let $n \ge 1$. Then, by the base case, by
subadditivity of the entropy, and by convexity and monotonicity of $\phi$:
$$
D^2(f) = \sum_{i=1}^n \E_x D^2\(f^{(x)}_i\)  =
4 \cdot \sum_{i=1}^n \E_x K^2\(f^{(x)}_i\)
\phi\(\frac{Ent\(\(f^{(x)}_i\)^2\)}{K^2\(f^{(x)}_i\)}\) =
$$
$$
4 K^2(f) \cdot \sum_{i=1}^n \E_x \frac{K^2\(f^{(x)}_i\)}{K^2(f)}
\phi\(\frac{Ent\(\(f^{(x)}_i\)^2\)}{K^2\(f^{(x)}_i\)}\) \ge
$$
$$
4 K^2(f) \cdot
\phi\( \frac{1}{K^2(f)} \cdot \sum_{i=1}^n \E_x Ent\(\(f^{(x)}_i\)^2\)\) \ge
4 K^2(f) \cdot \phi\(\frac{Ent\(f^2\)}{K^2(f)}\)
$$
\eprf

We proceed to derive (\ref{ineq:l_sob}) from (\ref{ineq:tech}). By homogeneity, we may assume $\E f^2 = 1$. Let $0 \le \rho
\le \log 2$. Consider the functional
$$
R[f] = \frac{D^2(f)}{Ent\(f^2\)}
$$
where $f$ ranges over the non-empty\footnote{Recall $Ent\(f^2\) \ge \E f^2 \log \frac{\E f^2}{\E^2 f}$ (\cite{FS}).} compact set of nonnegative non-zero functions satisfying $\E f^2 = 1$ and $Ent\(f^2\) \ge \rho n$, and the associated minimum
$$
m(\rho) = \min_f R[f]
$$
To complete the proof of (\ref{ineq:l_sob}) and of Theorem~\ref{thm:log-Sobolev}, we will show
\beqn
\label{ineq:mc}
m(\rho) \ge c(\rho),
\eeqn
where $c$ is the function defined in the third part of Lemma~\ref{lem:functions}.

Indeed, fix $\rho$ and let $m = m(\rho)$. Let $f$ be a function at which $R[f]$ attains its minimum, that is: $\E f^2 = 1$, $Ent\(f^2\) \ge \rho n$ and $D^2(f) = m Ent\(f^2\)$. Then
$$
D^2(f) = m Ent\(f^2\) \ge m \rho n.
$$
This means $K^2(f) = n \E f^2 - 1/4 D^2(f) \le \( 1 - (m\rho)/4\) n$, and therefore $\frac{Ent\(f^2\)}{K^2(f)} \ge \frac{4\rho}{4 - m\rho}$.

Recall $\phi$ is an increasing convex function on $[0,2 \log 2]$ with
$\phi(0) = 0$. Therefore the function $\tau(y) = \phi(y)/y$ is
increasing. This implies
$$
\phi\(\frac{Ent\(f^2\)}{K^2(f)}\) = \frac{Ent\(f^2\)}{K^2(f)} \tau\(\frac{Ent\(f^2\)}{K^2(f)}\) 
\ge \frac{Ent\(f^2\)}{K^2(f)}
\tau\(\frac{4\rho}{4 - m\rho}\)
$$
and, by Proposition~\ref{thm:technical},
$$
D^2(f) \ge 4 K^2(f) \phi\(\frac{Ent\(f^2\)}{K^2(f)}\) \ge 4
\tau\(\frac{4\rho}{4 - m\rho}\) Ent\(f^2\)
$$
which means
$$
m \ge 4 \tau\(\frac{4\rho}{4 - m\rho}\)
$$
The rest is simple algebra. Recall $\tau = \phi(y)/y$ and $\psi
= \phi^{-1}$. Since $\psi$ is increasing, the last inequality is equivalent to
$$
\psi\(\frac{m\rho}{4 - m\rho}\) \ge \frac{4\rho}{4 - m\rho}.
$$
Substituting  $y = \frac{m\rho}{4 - m\rho}$, this translates to
$\psi(y)/(y+1) \ge \rho$. Recall $\alpha = \(\frac{\psi(y)}{1 +
  y}\)^{-1}$. Since $\alpha$ is increasing, we obtain
$$
\frac{m\rho}{4 - m\rho} \ge \alpha(\rho)
$$
This is the same as
$$
m \ge \frac{4 \alpha(\rho)}{\rho(\alpha(\rho) + 1)}.
$$
This is equivalent to (\ref{ineq:mc}), completing the proof of (\ref{ineq:l_sob}) and of Theorem~\ref{thm:log-Sobolev}.
\eprf

\section{Proof of Lemma~\ref{lem:functions}}
Let
\beqn
h(t) = \frac12 (1-t)^2 \log (1-t)^2 + \frac12 (1+t)^2 \log (1+t)^2 - \(1+t^2\)\log\(1+t^2\).
\label{def:h}
\eeqn
In other words, $h(t) = \psi\(t^2\)$. We will start with some useful properties of the function $h$.

\lem
\label{lem:h-prop}
\begin{enumerate}
\item
$$
h' \ge h \ge 0
$$
\item
$$
\(1-t^2\) h' \ge t h''
$$
\end{enumerate}
\elem
\prf
We have $h'(t) = 2\cdot\((1+t)\log(1+t) - (1-t)\log(1-t) - t \log\(1+t^2\)\)$ and $h''(t) = \frac{4}{1+t^2} - 2\log\(\frac{1+t^2}{1-t^2}\)$.

The first claim of the lemma is easy. Nonnegativity of $h$ follows from nonnegativity of $\psi$, and
$$
h'(t) - h(t) = \(1-t^2\) \log(1+t) - (1-t)(3-t) \log(1-t) + \(1-t^2\) \log\(1+t^2\) \ge 0
$$
The second claim is somewhat harder. We have, rearranging and simplifying:
$$
\frac12 \cdot \(\(1-t^2\) h' - t h''\) = t^3 \log\(\frac{1+t^2}{1-t^2}\) + \(1-t^2\)\log\(\frac{1+t}{1-t}\) - \frac{2t}{1+t^2}
$$
That is, we need to show
$$
\(1+t^2\) t^3 \log\(\frac{1+t^2}{1-t^2}\) + \(1-t^4\) \log\(\frac{1+t}{1-t}\) \ge 2t
$$
In fact, an even stronger inequality
$$
t^3 \log\(\frac{1+t^2}{1-t^2}\) + \(1-t^4\) \log\(\frac{1+t}{1-t}\) \ge 2t
$$
is valid for $t \in [0,1]$. This is easy to check for $t=1$. To see this for $0 \le t < 1$, recall that for $-1 < x < 1$ holds $\log\(\frac{1+x}{1-x}\) = 2\sum_{k=0}^{\infty} \frac{x^{2k+1}}{2k+1}$. Substituting these series in the inequality above, we need to show
$$
t^3 \cdot \sum_{k=0}^{\infty} \frac{t^{4k+2}}{2k+1} + \(1-t^4\) \cdot \sum_{k=0}^{\infty} \frac{t^{2k+1}}{2k+1} \ge t,
$$
and this is easily verified by observing that all the higher coefficients of the power series on the left hand side are nonnegative.
\eprf

We pass to the proof of Lemma~\ref{lem:functions}

\noindent {\bf Claim 1}

\noindent

We will prove $\psi$ is concave by showing $\psi''$ is negative on $(0,1)$. We have for $0 < t < 1$:
$$
\psi''(t) = -\frac{1}{4t^{3/2}} \cdot h'\(\sqrt t\) + \frac{1}{4t}
\cdot h''\(\sqrt t\) < 0
$$
The last inequality follows from the second claim of Lemma~\ref{lem:h-prop}.

To see that $\psi$ is increasing, it suffices to verify $\psi' \ge 0$ at
$1$. Indeed,
$
\psi'(1) = \frac{1}{2} h'(1) = \log 2 > 0
$,
completing the proof of Claim 1.

\noindent {\bf Claim 2}

\noindent Let $\xi(t) = \frac{\psi}{1+t}$. We will verify $\xi' > 0$ on $(0,1)$. Indeed,
$$
\xi' = \frac{(1+t)\psi' - \psi}{(1+t)^2}
$$
so we need to check that $\psi' > \frac{\psi}{1+t}$. Substituting
$\psi(t) = h(\sqrt t)$, this amounts to checking
$$
h'(\sqrt t) \ge \frac{2\sqrt t}{1 + t} \cdot h(\sqrt t)
$$
Since $\frac{2\sqrt t}{1 + t} < 1$ on $(0,1)$, it suffices to
show $h' \ge h$, which is true by the first claim of Lemma~\ref{lem:h-prop}.

To show concavity of $\xi$, we will prove that
\beqn
-\xi'' > 2\xi,
\label{about-xi}
\eeqn
implying $\xi'' < 0$ on $(0,1)$. Direct calculation gives
$$
\xi''(t) = \frac{(1+t)^2 \psi'' - 2\((1+t)\psi' - \psi\)}{(1+t)^3}
$$
Simplifying, $-\xi'' > 2 \xi'$ reduces to
$$
2t \psi > 2t(1+t) \psi' + (1+t)^2 \psi''
$$
Since $\psi'' < 0$ we have $(1+t)^2 \psi'' < 4t \psi''$. Therefore, it suffices to prove
$$
\psi \ge (1+t) \psi' + 2 \psi''
$$
Again, since $\psi$ is concave with $\psi(0) = 0$, we have $\psi \ge t \psi'$. Therefore, we only need to prove
$$
-2 \psi'' \ge \psi'
$$
Writing this in terms of $h$, this is equivalent to $\(1-s^2\) h' \ge s h''$ which is given by the second claim of Lemma~\ref{lem:h-prop}.

\noindent {\bf Claim 4}

We will show this claim before the third claim of the lemma. That claim is somewhat more involved, and its proof is relegated to the end of this section.

Here we need to verify $$\frac{\alpha(t)}{1+\alpha(t)} = \frac12 - \sqrt{H^{-1}\(\log 2 - t\)\(1-H^{-1}\(\log 2 -t\)\)}$$ for all $t \in [0,\log 2]$. This is equivalent to
$$
H^{-1}\(\log 2 - t\) = \frac12 - \sqrt{\frac{\alpha}{1+\alpha}\(1 - \frac{\alpha}{1+\alpha}\)} = \frac{\(1-\sqrt{\alpha}\)^2}{2(1+\alpha)}
$$
Recall that $\alpha = \alpha(t)$ is defined to satisfy
$t = \frac{\psi(\alpha)}{1 + \alpha}$, where $\psi(\alpha) =
Ent\(f^2\)$, and $f$ is a function on $\{0,1\}$ with $g(0) = 1
- \sqrt \alpha$, $g(1) = 1 + \sqrt \alpha$. It is not hard to verify the identity
$$
\frac{\psi(\alpha)}{1 + \alpha} = \ln 2 - H\(\frac{\(1-\sqrt \alpha\)^2}{2(1 + \alpha)}\)
$$
for all $\alpha \in [0,1]$, and we are done.

\noindent {\bf Claim 3}

First, we show that $c$ is increasing. Direct computation gives that $c'$ is positive on $(0,\log 2)$ iff $t \alpha' > \alpha + \alpha^2$ on this interval. Both sides of this inequality are $0$ at
$0$, and we compare derivatives, that is, show  $t \alpha'' > 2 \alpha \alpha'$.

Since $\alpha$ is convex with $\alpha(0) = 0$, we have $\alpha \le t \alpha'$ in the interval.
Hence, it suffices to show $\alpha'' > 2\(\alpha'\)^2$.

Recall $\alpha = \xi^{-1}$. Consequently, $\alpha'(\xi(t)) =
\frac{1}{\xi'(t)}$, and $\alpha''(\xi(t)) =
-\frac{\xi''(t)}{\(\xi'(t)\)^3}$. Therefore, $\alpha'' >
2\(\alpha'\)^2$ is equivalent to $-\xi'' > 2\xi'$, which is given by (\ref{about-xi}).

It remains to show that $c$ is convex. This turns out to be significantly harder than the other proofs in this Section. We provide a somewhat sketchy argument below.

Direct computation shows that $c'' >0$ on $(0,\log 2)$ iff
\beqn
t^2 (1 + \alpha) \alpha'' + 2\alpha (1+\alpha)^2 >
2t^2\(\alpha'\)^2 + 2t (1+\alpha)\alpha'
\eeqn
First, we rewrite this inequality in terms of $\xi = \alpha^{-1}$. Let $t = \xi(x)$, that is $\alpha(t) = x$, $\alpha'(t) = \frac{1}{\xi'(x)}$, $\alpha''(t) = -\frac{\xi''(x)}{\(\xi'(x)\)^3}$. Substituting and simplifying, one gets
$$
-(1+x) \xi^2 \xi'' + 2x(1+x)^2 \(\xi'\)^3 > 2\xi^2 \xi' + 2(1+x) \xi \(\xi'\)^2,
$$
which has to hold for all $x$ in $(0,1)$.

Next, we rewrite this in terms of $\psi = (1+x) \xi$, obtaining
$$
(1+x)\psi^2\(-\psi''\) > 2\((1+x)\psi' - \psi\)^2 \(\psi- x\psi'\)
$$
Note that all the expressions in the brackets are positive, since $\psi$ is concave and $\xi$ is increasing, as we saw in the proofs of Claims 1 and 2 above. We simplify this inequality, replacing $\psi$ with $x \psi'$ on the left hand side and in the first term on the right hand side, and arriving to the stronger inequality
$$
x(1+x)\psi\(-\psi''\) > 2 \psi' \(\psi - x \psi'\)
$$
We rewrite this in terms of the function $h$, defined in (\ref{def:h}) above. As in the proof of Claim 1, expressing $\psi$ and its derivatives in terms of $h$, leads to the following equivalent inequality:
\beqn
2x(h')^2 > \(3-x^2\)hh' + x\(1+x^2\) h h''
\label{interms:h}
\eeqn
From now on we concentrate on the proof of (\ref{interms:h}). It will be convenient to write $h$ and its derivatives in terms of two new functions $L_1(x) = \log{\frac{1+x}{1-x}}$ and $L_2(x) = \log{\frac{1+x^2}{1-x^2}}$. Recalling the expressions for $h$ and its derivatives (as in the proof of Lemma~\ref{lem:h-prop} above, we have
\begin{itemize}
\item
$
h(x) = 2xL_1 - \(1+x^2\)L_2
$
\item
$
h'(x) = 2L_1 - 2xL_2
$
\item
$
h''(x) = \frac{4}{1+x^2} - 2L_2
$
\end{itemize}

Rewriting (\ref{interms:h}) in terms of $L_1$ and $L_2$, and simplifying, one arrives to
$$
\(3-x^2\)\(1+x^2\)L_1 L_2 + 2x\(1+x^2\)L_2 > 2x\(1-x^2\)L^2_1 + 4xL^2_2 + 4x^2L_1
$$
We expand both sides of this inequality as power series for $x \in (0,1)$. Recall that $L_1(x) = 2\sum_{k=0}^{\infty} \frac{1+x^{2k+1}}{2k+1}$, and, consequently, $L_2(x) = 2\sum_{k=0}^{\infty} \frac{1+x^{4k+2}}{2k+1}$. Therefore, both sides of this inequality have only odd terms.

Let the left hand side be equal to
$$
F(x) = 4 \cdot \sum_{k=0}^{\infty} \ell_{2k+1} x^{2k+1}
$$
and the right hand side be equal to $$
G(x) = 4 \cdot \sum_{k=0}^{\infty} r_{2k+1} x^{2k+1}
$$
We will argue that
\begin{enumerate}
\item
All the coefficients $\ell_{2k+1}$ and $r_{2k+1}$ are nonnegative.
\item
$\ell_1 = r_1 = 0$, $\ell_3 = r_3 = \ell_5 = r_5 = 4$.
\item
For all odd $k$ starting from $k = 3$:
$$
\ell_{2k+1} > r_{2k+1}~~~and~~~\ell_{2k+1} + \ell_{2k+3} > r_{2k+1} + r_{2k+3}
$$
\end{enumerate}

\noindent This will imply
$$
F(x) - G(x) = 4\cdot \sum_{k=3}^{\infty} \(\ell_{2k+1} - r_{2k+1}\) x^{2k+1} = 4\cdot \sum_{odd~k \ge 3} \(\(\ell_{2k+1} - r_{2k+2}\) - \(\ell_{2k+3} - r_{2k+3}\)x^2 \) \cdot x^{2k+1} >
$$
$$
4\cdot \sum_{odd~k \ge 3} \(\ell_{2k+1} - r_{2k+2}\)\(1-x^2\) \cdot^{2k+1} > 0,
$$
completing the proof of (\ref{interms:h}) and of Claim 3. Hence it remains to prove the properties of the coefficients.

In fact, the coefficients can be computed explicitly, which makes it possible to verify the required properties. We omit the (easy but cumbersome) details. For completeness sake, we do list explicit expressions for the coefficients below\footnote{Our apologies to the reader.}.
\begin{itemize}
\item
For an odd $k \ge 3$:
$$
\ell_{2k+1} = \(\frac{8k - 20}{(2k-3)(2k+1)}\) \cdot \sum_{m=1}^{(k-1)/2} \frac{1}{4m-3} ~~~+~~~ \frac{4}{2k-1} \cdot \sum_{m = 1}^{k-2} \frac{1}{2m+1} ~~~+~~~
$$
$$
\(\frac{3}{2k+1} +
\frac{2}{2k-1} - \frac{1}{2k-3}\) \cdot \sum_{m = 1}^{(k-1)/2} \frac{1}{2m-1} ~~~+~~~ \(\frac{1}{k} + \frac{3}{k(2k+1)} + \frac{6}{(2k-1)(2k+1)}\)
$$
and
$$
r_{2k+1} = \frac{2k+2}{k(2k-1)} ~~~-~~~ \frac{2}{k(k-1)} \cdot \sum_{m = 1}^{k-1} \frac{1}{2m-1}
$$
\item
For an even $k \ge 4$:
$$
\ell_{2k+1} = \(\frac{8k - 20}{(2k-3)(2k+1)}\) \cdot  \sum_{m = 1}^{k-2} \frac{1}{2m+1} ~~~+~~~ \frac{4}{2k-1} \cdot  \sum_{m=1}^{k/2} \frac{1}{4m-3}~~+~~
$$
$$
\(\frac{3}{2k+1} + \frac{2}{2k-1} - \frac{1}{2k-3}\) \cdot \sum_{m = 1}^{(k-2)/2} \frac{1}{2m-1} ~+~ \(\frac{1}{k-1} + \frac{6}{(2k-1)(2k+1)} +  \frac{10k-1}{(k-1)(2k-1)(2k+1)} \)
$$
and
$$
r_{2k+1} = \frac{8}{k} \cdot \sum_{m=1}^{k/2} \frac{1}{2m - 1} ~~~+~~~ \frac{2k+2}{k(2k-1)} ~~~-~~~ \frac{2}{k(k-1)} \cdot \sum_{m = 1}^{k-1} \frac{1}{2m-1}
$$
\end{itemize}

It remains to compute the coefficients for $k = 1, 2$. This is easily done directly, verifying the property 2 above.


\begin{thebibliography}{99}

\bibitem{ACKL}
A. Ashikhmin, G. Cohen, M. Krivelevich, S. Litsyn, {\sl Bounds on distance distributions in codes of known size},  IEEE Trans. Inform. Theory, vol. IT-51 2005, 250-258.

\bibitem{BJ}
A. Barg, D. B. Jaffe, {\sl Numerical results on the asymptotic rate of binary codes}, in {\bf Codes and Association schemes}, (A. Barg and S. Litsyn, eds.), Amer. Math. Soc., Providence 2001.

\bibitem{Beckner}
W. Beckner, {\sl Inequalities in Fourier Analysis}, Annals of Math.,
102(1975), pp. 159-182.

\bibitem{Bonami}
A. Bonami, {\sl Etude des coefficients Fourier des fonctiones de
  $L^p(G)$}, Ann. Inst. Fourier (Grenoble) 20:2 (1970), pp. 335-402.

\bibitem{Chavel}
I. Chavel, {\bf Isoperimetric inequalities}, Cambridge University
Press, 2001.

\bibitem{dels}
P. Delsarte, {\sl An algebraic approach to association schemes of coding theory}, Philips Res. Rep. Suppl. 10, 1973.


\bibitem{FS}
D. Falik, A. Samorodnitsky, {\sl A combinatorial proof of a theorem of
  Kahn, Kalai, and Linial}, CPC, to appear.

\bibitem{FT}
J. Friedman and J-P. Tillich, {\sl Generalized Alon-Boppana Theorems
and Error-Correcting Codes}, preprint, 2002.

\bibitem{Gr}
L. Gross, {\sl Logarithmic Sobolev inequalities}, Amer. J. of Math.,
97, 1975, pp. 1061-1083.

\bibitem{Harp}
L. H. Harper, {\sl Optimal numberings and isoperimetric problems on graphs}, J. Combin. Theory, 1, 1966, pp. 385-393.

\bibitem{Harper}
L. H. Harper, {\sl Optimal assignment of numbers to vertices},
J. Soc. Ind. Appl. Math., 12, 1964, pp. 131-135.

\bibitem{Hart}
S. Hart, {\sl A note on the edges of the $n$-cube}, Discr. Math., 14,
1976, pp. 157-163.

\bibitem{KL}
G. Kalai, N. Linial, {\sl On the distance distribution of codes}, IEEE Trans. Inform. Theory, vol. IT-41 1995, 1467-1472 (see, in particular, the extended version of the paper on http://www.cs.huji.ac.il/$\sim$nati/).

\bibitem{Ledoux}
M. Ledoux, {\sl Concentration of measure and logarithmic Sobolev
  inequalities}, Seminaire de Probabilites, XXXIII, Lecture Notes in
Math., 1709, Springer 1990, pp. 120-216.

\bibitem{lev}
V. I. Levenshtein, {\sl Krawtchouk polynomials and universal bounds for codes
and designs in Hamming spaces}, IEEE Trans. Inform. Theory, vol. IT-41,
1995, 1303-1321.

\bibitem{vLint}
J.H. van Lint, {\bf Introduction to Coding Theory}, Third edition.
Springer-Verlag, Berlin, 1999.

\bibitem{MRRW}
R. J. McEliece, E. R. Rodemich, H. Rumsey, Jr., and L. R. Welch,
{\sl New upper bounds on the rate of a code via the Delsarte-MacWilliams
inequalities}, IEEE Trans. Inform. Theory, vol. IT-23, 1977, 157-166.

\bibitem{NS}
M. Navon, A. Samorodnitsky, {\sl On Delsarte's linear programming
  bounds for binary codes}, Proceedings of FOCS 46.

\bibitem{NS1}
M. Navon, A.Samorodnitsky, {\sl Linear programming bounds via a covering argument}, Disc. and Comp. Geometry, to appear.

\bibitem{Szego}
G. Szeg\"o, {\bf Orthogonal Polynomials}, American Mathematical
Society, 1939.



\end{thebibliography}
\end{document}